\documentclass[graybox]{svmult}
\usepackage[T2A,OT4,T1]{fontenc}
\usepackage[utf8]{inputenc}
\usepackage[polish,russian,english]{babel}
\usepackage{mathptmx}       
\usepackage{helvet}         
\usepackage{courier}        
\usepackage{makeidx}         
\usepackage{graphicx}        
\usepackage{multicol}        
\usepackage[bottom]{footmisc}
\usepackage{amsmath,amssymb,amsfonts}        
\usepackage{bm}
\usepackage{color}
\usepackage{comment}

\RequirePackage[numbers]{natbib}%
\usepackage[]{hyperref}
\hypersetup{
    pdftitle={Holistic Decision-Making in Stopping Problems: Emphasizing Psychological Aspects},
    pdfauthor={Joanna Rymaszewska and Georgy Sofronov and Krzysztof Szajowski}, 
    colorlinks=true,
    urlcolor=blue,
    filecolor=magenta,
    citecolor=green, 
    linkbordercolor={1 1 1}, 
    citebordercolor={1 1 1},  
    urlbordercolor={ 1 1 1}  
  }
\RequirePackage[hyperpageref]{backref} 
    \renewcommand*{\backref}[1]{}  
    \renewcommand*{\backrefalt}[4]{
       \ifcase #1 
          No cited.
       \or
          Cited on p. #2.
       \else
          Cited on pp. #2.
       \fi}  	

\newcommand*{\MR}[1]{\href{http://www.ams.org/mathscinet-getitem?mr=#1&return=pdf}{MR #1}}
\newcommand*{\ZBL}[1]{\href{http://www.zentralblatt-math.org/zmath/en/advanced/?q=an:#1&format=complete}{Zbl #1}}

\newcommand*{\PMID}[1]{\href{https://pubmed.ncbi.nlm.nih.gov/#1/}{PMID: #1 }} 
 
\DeclareFontFamily{U}{mathx}{\hyphenchar\font45}
\DeclareFontShape{U}{mathx}{m}{n}{
      <5> <6> <7> <8> <9> <10>
      <10.95> <12> <14.4> <17.28> <20.74> <24.88>
      mathx10
      }{}
\DeclareSymbolFont{mathx}{U}{mathx}{m}{n}
\def\bbA\mathbf{A}

\def\bP{\bm{P}}
\def\bE{{\bm E}}

\def\bT{\mbox{\bf T}}
\def\bQ{\mbox{\bf Q}}

			\def\BbbE{\mathbb E}

\def\bbA{\mathbb A}

    \def\cB{{\mathcal B}}
		
	\def\cS{{\mathcal S}}

        \def\vtilde{\tilde v}

\def\whatf0{\widehat{f}^0}		
\def\cB{{\mathcal B}} 
\def\cF{{\mathcal F}}

    \def\DM{\textbf{DM}}

\def\MDM{\textbf{MDM}}

\def\bP{{\bf P}}
\def\bE{{\bf E}}

\def\bE{\bm{E}}
\def\bP{\bm{P}}
\def\bQ{\bm{Q}}
\def\bT{\bm{T}}

\def\BbbE{\mathbb E}
\def\BbbN{\mathbb N}
\def\cB{{\mathcal B}}
\def\cF{\mathcal{F}}
\def\cS{\mathcal S}

\newcommand*{\doi}[1]{\href{http://dx.doi.org/#1}{doi: #1}}
\makeindex             

\begin{document}
\title*{Holistic Decision-Making in Stopping Problems:\\ Emphasizing Psychological Aspects\thanks{Talk was presented on the MATRIX Research Program: Probabilistic Models in Evolutionary Biology and Game Theory (6-10 January 2025)
}}
\author{Georgy Sofronov, Joanna Rymaszewska, and Krzysztof J. Szajowski}
\institute{Georgy Sofronov \at  School of Mathematical and Physical Sciences, Faculty of Science and Engineering, Macquarie University NSW 2109 Australia, \email{georgy.sofronov@mq.edu.au} \and  Joanna Rymaszewska \at Faculty of Medicine, WUST, Hoene-Wrońskiego 13c, 50-372 Wrocław, Poland, \email{joanna.rymaszewska@pwr.edu.pl}\and Krzysztof J. Szajowski \at Faculty of Pure and Applied Mathematics, Wrocław University of Science and Technology (WUST), Wybrzeże Wyspiańskiego 27, 50-370 Wrocław, Poland, \email{Krzysztof.Szajowski@pwr.edu.pl}}
%
%
\maketitle
\vspace{-5ex}
\abstract{Our research is closely related to \emph{ontological studies} in mathematics. It provides crucial insights into the nature of decisions and strategies characterized by Markov moments.\\
In a stopping game, a holistic decision-maker would \emph{evaluate comprehensive information} by assessing the probabilities of various outcomes and their associated payoffs. This involves understanding the current state, historical data, and potential future scenarios. Such a decision-maker must also consider \emph{strategic interactions} by anticipating and accounting for the strategies of other players. They must be flexible in adapting their strategy as the game evolves and able to integrate uncertainty by incorporating risk preferences and tolerances. They would perform scenario analysis to evaluate the impact of different stopping times under varying conditions.\\
The goal of this modeling and its implementation in psychological practice is to introduce a novel method for assessing the state of players, leveraging deviations from rational strategies as diagnostic indicators of their psychological and decision-making profiles. The details of other models will be subject to contributed papers. The article presents the theoretical basis for combining various factors when modeling decision-making processes. The original title is "Rationality, Deviation, and Diagnosis: A Holistic Approach to Stopping Games" and will be used when it is possible to describe and interpret the results of the experiments we write about in the last section of the paper.
}

\section{\label{KSzSec:1}Strategic Interaction, and Uncertainty in stopping games}
    A holistic decision-maker (\DM) considers the entire system and all relevant factors when making a decision, rather than focusing on isolated parts. 
    In a stopping game, this means evaluating the potential results of stopping the game at different points (possible outcomes), understanding how different elements of the game influence each other (interdependencies), and looking beyond immediate gains or losses to consider future implications (long-term consequences).
    
    Participants in the decision-making process would \emph{evaluate comprehensive information} to assess the probabilities of various outcomes and their associated rewards and understand the current state, historical data, and potential future states.
    They also should \emph{consider strategic interactions}, anticipate and account for the strategies of other players. 
    Such players should be flexible in adjusting their strategy based on how the game evolves and incorporate risk preferences and tolerances (\emph{integrate uncertainty}). In addition, they would perform scenario analysis to understand the impact of different stopping times under various conditions.
    
    A complete approach to modeling often means focusing on the basic parts of math concepts. This places our research in the exciting area of studying the nature and structure of mathematics (v. \emph{Mathematics and reality} by~\citeauthor{Len2010:Reality}~(\citeyear{Len2010:Reality})).

A prototypical illustration is experimental psychology research, which can be regarded as ontological in the context of the secretary problem (v.~\citeauthor{BearMurRap2005:SP}~(\citeyear{BearMurRap2005:SP}), \citeauthor{ErnstSza2021:Average}~(\citeyear{ErnstSza2021:Average})) due to its capacity to furnish invaluable insights into the essence of decisions and strategies that exhibit Markov characteristics.

Consider two decision-makers, \textbf{A} and \textbf{B}, in a collaborative setting. Decision-maker \textbf{A} might aim to maximize short-term gains, while decision-maker \textbf{B} may aim for long-term stability. They must decide when to stop observing a fluctuating market. They agree to stop if either
\begin{itemize}
\item both agree that the market conditions are optimal for stopping, or
\item only one of the players identifies a critical condition (e.g., market crash), which triggers an immediate stop.
\end{itemize}
Defining clear stopping rules for each participant and understanding how their decisions impact the overall process are crucial for formulating an effective multiple-stop strategy.

\subsection{Detailed aspects of Multilateral Decision-Making under uncertainty} In an increasingly interconnected world, \emph{multilateral decision-making} (\MDM) has become a crucial process for addressing complex problems that involve multiple stakeholders. Whether in international diplomacy, economic policy, or large-scale organizational strategy, decisions are rarely made in isolation. However, {\MDM} is often conducted under conditions of uncertainty, where unpredictable external factors, incomplete information, and diverse stakeholder interests complicate rational decision-making.

To enhance decision-making effectiveness, it is essential to integrate mathematical models with a deep understanding of human behavioral factors. This approach ensures that decisions remain rational, adaptive, and robust in dynamic environments. The following discussion explores the sources of uncertainty in \MDM, the mathematical techniques available for managing it, and the role of behavioral considerations in improving decision outcomes.

\subsubsection{Understanding Uncertainty in \MDM} Uncertainty in decision-making arises from multiple sources, making it a fundamental challenge in \MDM. Many real-world phenomena, such as economic fluctuations, climate patterns, and resource availability, are inherently uncertain due to natural variability. Additionally, decision-makers often lack comprehensive data about the systems they are dealing with, leading to gaps in understanding. For example, in geopolitical negotiations, incomplete intelligence about another nation's strategic intentions can significantly impact decision outcomes. Beyond these factors, unforeseen events such as financial crises, technological breakthroughs, or sudden policy shifts can render prior decisions obsolete or suboptimal. These uncertainties require continuous monitoring and flexibility in decision-making frameworks.

A structured approach to identifying and categorizing these sources allows decision-makers to develop strategies to mitigate their impact. While uncertainty cannot always be eliminated, understanding its nature enables more effective planning and response strategies.

\subsubsection{Mathematical Approaches to Managing Uncertainty} To enhance rational decision-making under uncertainty, various quantitative methods provide a systematic way to assess risks, evaluate alternatives, and anticipate potential outcomes. Probabilistic models, for example, assign probabilities to different outcomes, allowing decision-makers to weigh their choices more effectively. Bayesian models offer a particularly valuable approach by updating probabilities as new data emerges, making them well-suited for dynamic decision-making environments.

Information gathering and sharing play a crucial role in managing uncertainty. Expanding data collection through surveys, real-time monitoring, and intelligence gathering enhances decision accuracy. Moreover, ensuring transparency and open communication among stakeholders reduces asymmetry and promotes collaboration. Scenario analysis and sensitivity testing further strengthen decision processes by simulating multiple possible futures and assessing how variations in input data affect outcomes. These tools collectively provide a rational foundation for decision-making, helping to mitigate uncertainty and support more informed choices.

While mathematical models offer structured decision-making frameworks, they must be complemented by an understanding of human behavioral factors, which often introduce additional complexities into multilateral decision processes.

\subsubsection{Behavioral Considerations in \MDM} Although mathematical models enhance objectivity, decision-making is ultimately a human-driven process influenced by psychological and social factors. Cognitive biases, for instance, frequently affect judgment, leading decision-makers to favor information that aligns with preexisting beliefs or to rely too heavily on initial pieces of data. Addressing these biases is essential to ensuring more rational and balanced decision-making.

In addition to cognitive limitations, risk preferences among stakeholders influence decision outcomes. Some decision-makers exhibit greater tolerance for uncertainty, while others prioritize stability. Because of these differences, {\MDM} frameworks should be adaptable to evolving goals and shifting circumstances. Psychological and social dynamics also shape decision-making processes. Stress, fatigue, and group interactions can influence judgment, sometimes leading to suboptimal choices. Furthermore, disruptive elements, such as conflicting stakeholder interests or adversarial actors, can destabilize decision processes. Effective mitigation strategies, including negotiation, alliance-building, and conflict resolution, help manage such disruptions and preserve decision coherence.

Integrating behavioral insights with mathematical models ensures that decision-making remains both analytically sound and practically viable. By considering the human dimension of decision-making, {\MDM} can become more resilient, inclusive, and responsive to real-world complexities.
\medskip

\noindent{\bf In summary}, multilateral decision-making under uncertainty is a complex process whose analysis — and consequently, modeling — requires both mathematical precision and behavioral awareness. By systematically categorizing uncertainty, applying probabilistic models and scenario analyses, and accounting for human cognitive and social factors, decision-makers can improve both the rationality and effectiveness of their choices.

Since uncertainty is an inherent feature of dynamic systems, continuous learning, adaptation, and feedback loops are essential for refining decision strategies over time. In an era of rapid change, those who can effectively balance data-driven insights with human adaptability will be best prepared to navigate uncertainty and make informed, strategic decisions. This highlights the importance of quantifying the adaptability of decision-makers and identifying the factors that influence the decision-making process. 

\emph{At the same time, observing decision-making behavior provides valuable insights into distortions in the correct reception and evaluation of stimuli. This underscores the significance of describing these distortions through behavioral observations in simulated scenarios.}

\subsection{Utilizing Mathematical Models} We start with modeling the phenomenon. This includes system representation, which involves creating accurate models that capture the system's dynamics, including variables, constraints, and interdependencies. Simulation techniques allow us to predict outcomes under different scenarios, helping anticipate potential problems and plan accordingly.

Next, we define objectives and measure impacts. This requires clearly defining objective functions that reflect the decision-makers' goals and can be optimized using mathematical techniques. We also develop key performance indicators (KPIs) that provide measurable metrics to assess the effectiveness of actions taken.

A critical limitation of most initial analyses is that they often ignore the rationality and acceptability of decision-makers. While mathematical models enhance decision-making, human judgment remains a crucial factor due to cognitive biases. Common biases, such as confirmation bias and anchoring, can distort decisions. To mitigate their effects, we must integrate decision support systems that incorporate model outputs and provide balanced, unbiased information.

Behavioral considerations are essential in decision-making. This includes aligning the process with the risk tolerance and preferences of decision-makers. Psychological factors such as stress, fatigue, and group dynamics can influence decisions, so it is important to account for them.

To ensure robustness, we conduct sensitivity and scenario analyses to assess how deviations in model assumptions affect outcomes. Finally, we implement a feedback loop, where decisions are continuously reviewed and adjusted based on new information and observed results.

 \section{Optimal sequential search}
 \subsection{The optimal stopping problems} 
 An optimal stopping problem with one stopping option can be presented as follows. Let $(X_n,{\cF}_n,\bP_x)_{n=0}^N $ be a homogeneous Markov chain defined on a probability space $(\Omega,{\cF},\bP)$ with a state space, $(\BbbE,\cB)$ and let $f:\BbbN\times\BbbE\rightarrow\Re$ be a sequence of real-valued $\cB$-measurable functions. Assume that horizon $N$ is finite.
 
\begin{definition}\label{MEKSzST}
The random variable $\tau:\Omega\rightarrow \BbbN$ such that $\{\omega:\tau(\omega)=n\}\in\cF_n$ for every $n\in \BbbN$ is called a Markov moment.  The Markov moment $\tau$ almost everywhere (\textit{a.e.}) finite is a stopping time. Let ${\cS^N}$ be the aggregation of Markov times. $\bP_x(\tau\leq N)<1$ for some $\tau\in{\cS^N}$. Denote $\cS_n^N=\{\tau\in\cS^N:\tau\geq n\}$. Let 
\begin{equation}\label{MEKSzOSV}
v(n,x)=\sum_{\tau\in\cS_n^N}\bE[f(\tau,X_{\tau})|X_n=x].
\end{equation}
The function $v(0,x)$ is the value of the optimal stopping problem. 
\end{definition}

Let us define  
\begin{align}\label{MEKSzOperT}
\bT f(n,x)&=\bE[f(n+1,X_{n+1})|X_n=x] \text{ (the mean operator),}\\
\bQ f(n,x)&=\max\{f(n,x),\bT f(n,x)\} \text{ (the maximum operator)}.\label{MEKSzOperQ}
\end{align}
\begin{theorem}
The value function $v(k,x)$ fulfills the equation
\begin{align*}
v(N,x)&=f(N,x),\\
v(n,x)&=\max\{f(n,x),Tv(n,x)\} \text{ for $n\in\BbbN$}.
\end{align*}
The Markov time $\tau^\star=\inf\{n\in\BbbN: X_n\in A_n\}$, where $A_n=\{x\in\cB:f(n,x)\geq \bT v(n,x)\}$, is optimal.
\end{theorem}

\subsection{Some discussion from the psychology-an ambivalence}\label{KSzPsy1}
The aim of the study is to analyze the avoidance--striving psychological conflict (v. \cite[Chapt. 13]{Sok2005:PsychDecRyz}, \cite{Koz1981:PDT}) in the selection process. It is the result of mentally uncomfortable ambivalence and cognitive dissonance that can lead to avoidance, procrastination, or deliberate attempts to resolve the 
{\it ambivalence}.
The idea of ambivalence was coined by Paul Eugen Bleuler in the beginning of the XX century when he researched on schizophrenia (v. \cite{FryKiej2008:Schizo}). It may or may not be perceived as psychologically unpleasant when the positive and negative aspects of the subject are simultaneously present in the person's mind (cf. \cite{NewbyClark2002:ThinkingAC}, \cite{SonEwo2014:Metacognitive}). In latest book \cite{Pin2021:Rationality} alludes to the relationship between this problem and rationality (Vide sections: \emph{Conflicts among Goals} and \emph{Conflicts among Time Frames} of \cite[Chapt. 2]{Pin2021:Rationality}).

\textbf{DM} tries to avoid risks and undesirable situations and, at the same time, strives to achieve the desired and assumed goals. Adopting this, the final effect of selection has two general states: one is equated with loss, danger, and worst results; the second is usually equated with the achievement of the assumed goal, which brings specific benefits. The choice is a strive-avoidance conflict. It expresses the split between "fear and greed" (v. \citeauthor{Coo1975:Portfolio}~(\citeyear{Coo1975:Portfolio}) introducing by \citeauthor{Mar1959:Portfolio}~(\citeyear{Mar1959:Portfolio}) to psychology, cf. also \citeauthor{CapKop2014:Portfolio}~(\citeyear{CapKop2014:Portfolio}), \citeauthor{BieRut20002:Credit}~(\citeyear{BieRut20002:Credit})).

Consider the model of sequential analysis of options in order to choose the best option that they are all different — so only one is the best and the choice of this one will satisfy him if all of them are not exhausted. The choice is risky because the options come in random order. However, it is easy to assess the probability that the analyzed object is the best. The results of the research on this model are known in the literature, both theoretical and experimental. In experimental research, attempts were made to assess to what extent the decision-makers correctly assess the chances that the tested object is the best. This was related to the expansion of information about the analyzed objects, which is justified by the real problems to which this type of model is applied, but does not answer the question about the role of motivational mechanisms, including personality-conditioned goals in the decision problem, in which the effect of the choice is to satisfy goals or failure. Here, how the way decisions are made is influenced by the conflict between "fear and hope". 

Such selection tasks can be found in the transplant profession. To put it simply, when an organ for transplantation appears, the pending patients are checked for compliance. One of the most important factors in the success of the procedure is choosing the best patient for the now available hospital materials. It is one of the most important factors. There is a risk of transplant organ rejection with no compatible patient selection. This is a situation we fear. As a result of such a "matching" procedure, if a compatible patient is not selected early enough, we will lose the opportunity to perform the procedure due to the lack of a candidate for the procedure—this is another circumstance we are concerned about. The candidate's material and the chance of a better one will be significantly lower than the chance that the analyzed patient is the best (and we do not know this at the time of the examination because the remaining patients are to be examined later (v. \citeauthor{Kras2010:organ}~ (\citeyear{Kras2010:organ})). 

\subsubsection{\label{MEKSzPen}Reward for right choice \& penalty otherwise}
The selection problem of the one decision maker when 
\begin{itemize}
\item the success is granted by $\alpha$,
\item the failure is penalized by $\beta$, and  
\item the case of no choice gives the penalty $\gamma$.
\end{itemize} 

The state space $\BbbE=\{1,2,\ldots,n\}\cup\{\partial\}$. When the \DM{}  reaches the state $\partial$, it is automatically stops with penalty $\gamma$. Let $f:\BbbE\rightarrow \Re$. The operator $\bT f(x)=\int_{\BbbE}f(t) p(x,dt)=\sum_{j=k+1}^np(x|t)f(t)+f(\partial)p(x|\partial)$. The option under consideration at moment $k$ which is the candidate for the winning choice should be the relative first. The process is called to be in state $k$, where $1\leq k\leq n$, if the decision maker face the $k$-th option being a candidate. Let us assume that no stop has been made until the state $k=s$. Choosing this option, the expected reward is equal to
\begin{equation}\label{reward1}
g(s)=\left\{\begin{array}{ll}
(\alpha+\beta)\frac{s}{n}-\beta &\text{ if $s\in\{1,2,\ldots,n\}$},\\
-\gamma &\text {otherwise.}
\end{array}
\right.
\end{equation}

Denote $V(k)$ the maximum expected payoff when the process is in the state $k$. The optimality principle of dynamic programming gives the equation:
\begin{align}\nonumber
V(k)&=\max\{g(k),\bT V(k)\}\\
&=\max\{(\alpha+\beta)\frac{k}{N}-\beta,\sum_{j=k+1}^N\frac{k}{j(j-1)}(\gamma+V(j))-\gamma\},\label{optrewardeq}
\end{align}
for $k=1,2,\ldots,n$, with $V(n)=\alpha$.

It is worth emphasizing that equation \eqref{optrewardeq} is not the same as the one in the classical secretary problem because of the different reward here. By this fact, the solution will be different, and the probability of selecting the best option will also be different. Let us define
\begin{equation}\label{Hkn}
H_{k,n}=\sum_{j=k}^{n-1}\frac{1}{j}.
\end{equation}
\begin{theorem}[The value and strategy of the problem]\label{OneDMNIBCP} The solution of \eqref{optrewardeq} is 
\begin{equation}\label{MEPV}
V(k)=\left\{
    \begin{array}{ll}
        \frac{k^\star-1}{N}\left((\alpha+\beta)H_{k^\star,N}+\beta\right)-\beta&\text{for $k=1,2,\ldots,k^\star-1$},\\
        (\alpha+\beta)\frac{k}{N}-\beta&\text{for $k=k^\star,\ldots,n$,}
    \end{array}
\right.
\end{equation}
where $k^\star$ is the integer $k$ for which $H_{k-1,N}\geq\frac{\alpha+\gamma}{\alpha+\beta}>H_{k,N}$ and $H_{k,N}$ is given by \eqref{Hkn}. 
\end{theorem}

\subsubsection{\label{MEKSzPen}Optimality and asymptotic behavior}
In other words, the optimal policy is to reject the first $k^\star-1$ options and then to accept the first candidate thereafter. 
This is the optimal stopping region, since the conditions for the monotone case of \citeauthor{ChoRobSig1971:Great}~(\citeyear{ChoRobSig1971:Great}) (see also \citeauthor{Shi1969:StatPosAnal}~(\citeyear{Shi1969:StatPosAnal})) are fulfilled. 
\begin{corollary}\label{MEKSzOSCorr}[Asymptotics] 
When $n\rightarrow\infty$ we have:
\begin{align}\nonumber
t^\star&=\lim_{n\rightarrow\infty} \frac{k^\star}{n}=\exp\left(-\frac{\alpha+\gamma}{\alpha+\beta}\right),\\
     V(0^{+})&=\lim_{n\rightarrow\infty} \frac{k^\star-1}{n}\left((\alpha+\beta)H_{k^\star,n}+\beta-\gamma\right)-\beta\nonumber\\
     &=(\alpha+\beta)\exp\left(-\frac{\alpha+\gamma}{\alpha+\beta}\right)-\beta,\nonumber\\
\label{KSzPwin}    {\bP}(\text{\DM{} win})&=\frac{\alpha+\gamma}{\alpha+\beta}\exp\left(-\frac{\alpha+\gamma}{\alpha+\beta}\right).
\end{align}
\end{corollary}
\subsubsection{Solution for various profile of the decision-maker} The model under consideration describes a situation in which the decision-maker distinguishes between more than one goal. The essential one is to indicate the global maximum, but it also takes into account the behavior in which, by deciding to wait longer for the right observation, we do not have a chance to choose. The details are formulated in Section~\ref{MEKSzPen}. The general description of the optimal solution in these problems can be done by subset of the state space.  Define $A^{(1)}_n=\{x\in\BbbE:g(n,x)=\bT\vtilde_{N}(n+1,x)\}$, and 
\begin{equation}\label{OFirstS}
\tau^{(1)}=\inf\{j\leq N: X_j\in A^{(1)}_j \}.
\end{equation}
We have the possibility of describing the profile \DM{} using various values of the parameters $\alpha$, $\beta$, and $\gamma$: 
\begin{description}
\item[ $\beta=\gamma=0$]  This is the classical secretary problem with $\tau_1^\star$ of threshold type. The threshold $k^\star\cong \lfloor ne^{-1}\rfloor +1$ (see \citeauthor{GilMos1966:maximum}~\cite{GilMos1966:maximum}).
\item[$\gamma=0$] It is the case without additional penalty for lack of selection. The stopping set consists of states $k$ such that $\frac{\alpha}{\alpha+\beta}>H_{k,n}$. 
\end{description}

\subsubsection{\label{MEKSzMTPenBCP}Duration of the decision process} The optimal strategy in this case is presented in Section~\ref{MEKSzPen} and is explicitly given in Theorem~\ref{OneDMNIBCP} and Corollary~\ref{MEKSzOSCorr} (v.~\citeauthor{ErnstSza2021:Average}~ (\citeyear{ErnstSza2021:Average})). The optimal strategy is a threshold one. For $\tau^\star$ given in the Corollary~\ref{MEKSzOSCorr}, we have $\bP(\omega:\tau^\star\leq N)=\frac{k^\star-1}{N}H_{k^\star,N}\cong \exp\{-\frac{\alpha+\gamma}{\alpha+\beta}\}<1$. To describe the decision-making time, we will modify this optimal Markov moment and create $ \tilde{\tau}^\star$ which are equal to the optimal Markov moments on the set $ \{\omega: \tau^\star<n \} $, and on the complement of this event they are equal to $n$. We can interpret $ \tilde{\tau}^\star $ as the time to stop on one state or reach the last observation at $n$. The full decision-making process ends at $ \tilde{\tau}^\star $ and we are interested in the mean of this random variable. Let us denote $m(k)=\bE(\tilde{\tau}^\star|\cF_k)$, $k\leq n$.  
\begin{align}\label{KSzDurMk}
m(k)&=\left\{
\begin{array}{ll}
\sum_{j=k+1}^N\frac{k}{j(j-1)}j+N(1-\sum_{j=k+1}^N\frac{k}{j(j-1)}) & \mbox{ if $ k>k^\star$,}\\
\sum_{j=k^\star+1}^N\frac{k^\star}{j(j-1)}j+N(1-\sum_{j=k^\star+1}^N\frac{k^\star}{j(j-1)})   & \mbox{ if $ k\leq k^\star$;}
\end{array}
\right.\\
\nonumber&=\left\{
\begin{array}{ll}
k H_{k,N}+k& \mbox{ if $ k>k^\star$,}\\
k^\star H_{k^\star,N}+k^\star   & \mbox{ if $ k\leq k^\star$.}
\end{array}
\right.
\end{align}

\begin{corollary}\label{MEKSzOptThers}When $n\rightarrow\infty$ we have:
\begin{align}
\nonumber\hat{m}(t)&=\lim_{N\rightarrow\infty} \frac{m(k)}{N}=\left\{
\begin{array}{ll}
-t\ln(t)+t& \mbox{ if $ t>t^\star$,}\\
-t^\star\ln(t^\star)+t^\star   & \mbox{ if $ t\leq t^\star$;}
\end{array}
\right.\\
\label{KSzAsDur}\hat{m}(0^+)&=\lim_{t\downarrow 0}\hat{m}(t) =\Big[1+\frac{\alpha+\gamma}{\alpha+\beta}\Big]\exp\left(-\frac{\alpha+\gamma}{\alpha+\beta}\right).
\end{align}
\end{corollary}
For $\alpha=1$ and $\beta=\gamma$ we have the result by \citeauthor{Yeo1997:Duration}~\cite{Yeo1997:Duration}.

Further details and numerical analysis can be found in \citeauthor{SofSza2025:Multiple}~(\citeyear{SofSza2025:Multiple}).    


\section{Experiments}
\subsection{Sequential choice of the global maximum}
The implementation of an experiment that meets the assumptions of selecting the best object requires the preservation of the assumptions of this problem. We assume that the observer keeps track of the appearing elements, which are a priori ordered according to a random permutation. The goal is to select the best element when it is shown to the observer. No additional knowledge at the time of presentation of the $k$-th observation is provided, except for the ability to determine its relative rank among those already presented. 

The number $N$ of all objects available in the experiment is also known, and it is known that no two objects are identical among them.  Therefore, we cannot present permutations of the number $1,2,\ldots,N$ with the stated goal of choosing the smallest number. Especially if the experiment is to be repeated many times by the same selector. Thus, it is necessary to create a set of objects that, when they are all known, can be easily ordered. However, in a sequential procedure, not all of them are known, but only those that are presented. As for the others, nothing but their number can be inferred. Even when the experiment is repeated many times. An algorithm that guarantees this can, for example, have a database of numbers up to $1$ to $5^N$ from which we randomly select a set of $N$ numbers before presentation. These numbers are presented sequentially, and the observer succeeds when he or she stops at the observation that is the maximum in this selected subset. To reduce the information as to the maximum possible value in the set of observations, we draw randomly from among $1$ to $a^N$, where $a\in \{2,3,5,7,11\}$ is randomized for each trial. 

\subsection{Collection of the post experimental data}
In the version of the decision experiment under consideration, there are two possible outcomes: the indication of a global maximum or no such indication. In addition, it is possible to distinguish the situation that the lack of success is due to the achievement of the maximum number of observations without indication or the indication of the wrong observation. Thus, after a series of experiments, there should be available:
\begin{itemize}
    \item number of experiments (cf. Lemma~\ref{MEKSzOSCorr});
    \item duration of each experiment (cf. formulae~\eqref{KSzDurMk}, \eqref{KSzAsDur});
    \item number of experiments successfully completed (cf.~\eqref{KSzPwin});
    \item number of experiments in which the observer reached the last observation and did not indicate a candidate
    \item the number of experiments in which the observer pointed to the wrong candidate. 
\end{itemize}

\section{Conclusion} We highlight experiments whose mathematical models are well-known and thoroughly studied. The significance of individual assumptions is well understood, meaning that the mathematical model can quantify deviations from these assumptions in an experiment. We argue that such deviations may be related to the way the decision-making problem is analyzed by the decision-maker. Therefore, the experimental results should allow for the detection and measurement of such regulatory disruptions. Studies such as those by \citeauthor{SteSeaRap2003:Heuristic}~(\citeyear{SteSeaRap2003:Heuristic}), \citeauthor{mizrLauf2024:NeuroSDM}~(\citeyear{mizrLauf2024:NeuroSDM}) also suggest similar cause-and-effect relationships. The proposed parameterization enables a more precise examination of these dependencies.

Analyzing decision-making models in group interactions can provide valuable insights into participants' mental health. Different decision-making approaches reveal how emotions and psychological well-being influence individual choices. The traditional rational choice model assumes that decisions are based on full information and objective analysis, but emotional challenges like anxiety or depression can impair judgment and rational decision-making. Game theory, including Nash equilibrium, helps assess one's ability to cooperate and anticipate others' actions—skills that may be compromised in individuals struggling with mental health issues. However, the usage of such models in diagnosis is limited due to multiple solutions (cf.~\citeauthor{RavSza1992:NonZero}~(\citeyear{RavSza1992:NonZero}), \citeauthor{RadSza90:Sequential}~(\citeyear{RadSza90:Sequential}), \citeauthor{NowSza1999:ISDG}~(\citeyear{NowSza1999:ISDG})). More profitable seems hierarchical decision-making models, such as Stackelberg equilibrium, which highlight how leadership styles impact the well-being of a group (cf.~\citeauthor{SkarSzaj2024:Stackelberg}~(\citeyear{SkarSzaj2024:Stackelberg}). Cognitive biases and heuristics play a major role in shaping decisions, often making individuals prone to errors, especially when their emotional state affects perception and reasoning. Social interactions and emotions further shape group dynamics, as stress, anxiety, or depression can significantly alter both collective and individual decision-making processes. Beyond mathematical models, qualitative approaches such as interviews and observations offer deeper insights into how individuals perceive their own decisions and the emotions driving them. Understanding these factors can help bridge the gap between decision science and mental health, leading to more effective support strategies for those facing psychological challenges.


\end{document}